\documentclass[12pt,a4paper]{article}
\RequirePackage{etex}
\usepackage{amsfonts, amsmath, amssymb,amsthm,graphics,mathrsfs,verbatim, CJK}
\usepackage{makeidx}
\usepackage{palatino}
\usepackage{titletoc}
\usepackage{tocloft}
\usepackage{titlesec}
\usepackage{stmaryrd,color}
\usepackage[colorlinks=true]{hyperref}
\usepackage{graphicx}
\usepackage{pgfplots}
\usepackage{tikz}
\usepackage{pstricks,pst-node,multido,ifthen,calc}

\usetikzlibrary{arrows}
\usetikzlibrary{calc}

\pgfplotsset{compat=1.16}

\titleformat{\section}{\Large}{}{0.5em}{}
\makeindex

\normalsize

\newtheorem{mthm}{Theorem}[section]
\newtheorem{mylem}[mthm]{Lemma}
\newtheorem{myprn}[mthm]{Proposition}
\newtheorem{mycor}[mthm]{Corollary}
\newtheorem{mydef}[mthm]{Definition}
\newtheorem{myrem}[mthm]{Remark}
\newtheorem{mycon}[mthm]{Construction}
\newtheorem{myeg} [mthm]{Example}
\newtheorem{myque} [mthm]{Question}
\newtheorem{myalg} [mthm]{Algorithm}
\newtheorem{myconj} [mthm]{Conjecture}
\newenvironment{thm}{\begin{mthm}}{\end{mthm}}

\newenvironment{prop}{\begin{myprn}}{\end{myprn}}

\def \Lemma #1 {\vs{2mm}\nin {\bf Lemma #1} }
\def \Prop #1 {\vs{2mm}\nin {\bf Proposition #1} }
\def \Th #1 {\vs{2mm}\nin {\bf Theorem #1} }
\def \Cor #1 {\vs{2mm}\nin {\bf Corollary #1} }
\def \Proof{\vs{2mm}\nin {\bf Proof.~}}
\def \Def #1 {\vs{2mm}\nin {\bf Definition #1} }
\def \part #1 {\hfil\break\hglue 12pt {\rm (#1)~}}

\def\D{\Delta}

\def\eset{\emptyset}

\def\fs{\footnotesize}

\def\Llra{\Longleftrightarrow}

\def\nin{\noindent}

\def \qed {~\vrule height6pt width 6pt depth 0pt}

\def\seq{\subseteq}

\def\vs{\vspace*}
\def\vi{\varphi}

\def\CM{ Cohen-Macaulay }
\def\iff{ if and only if }
\def\tfae{ the following statements are equivalent}
\def\vd{vertex-decomposable}

\makeatletter
\newcommand{\eqnum}{\refstepcounter{equation}\textup{\tagform@{\theequation}}}
\makeatother

\begin{document}
\title{
\bf\LARGE  The Cohen-Macaulay Property of $f$-ideals\thanks{This research was supported by the Natural Science
Foundation of Shanghai (No. 19ZR1424100), and the National Natural Science Foundation of China (No. 11971338) } }

\author{ A-Ming Liu\thanks{  aming8809@163.com},\,\,\,Jin Guo\thanks{guojinecho@163.com}\\
{\fs School of Science, Hainan University}\\
Tongsuo Wu\thanks{Corresponding author. {\small
tswu@sjtu.edu.cn}}\\
{\fs School of Mathematical Sciences, Shanghai Jiao Tong University}\\}

\date{}
\maketitle

\begin{center}
\begin{minipage}{12cm}

\vs{3mm}\nin{\small\bf Abstract} {\fs For positive integers $d<n$, let $[n]_d=\{A\in 2^{[n]}\mid |A|=d\}$ where  $[n]=:\{1,2,\ldots, n\}$. For a pure $f$-simplicial complex $\D$ such that ${\rm dim}(\D)={\rm dim}(\D^c)$ and $\mathcal{F}(\D)\cap \mathcal{F}(\D^c)=\eset$, we prove that the facet ideal $I(\D)$ is Cohen-Macaulay if and only if it has linear resolution. For a $d$-dimensional  pure $f$-simplicial complex $\D$ such that $\D'=:\langle F\mid F\in [n]_d\smallsetminus \mathcal F(\D)\rangle$ is an $f$-simplicial complex, we prove that $I(\D^c)$ is Cohen-Macaulay if and only if $I(\D')$ has linear resolution.}

\vs{3mm}\nin {\small\bf Key Words} {\small $f$-ideal; Cohen-Macaulay; Newton complement dual; linear resolution}

\vs{3mm}\nin {\small 2010 AMS Classification:} {\small Primary: 13H10; 05E45; Secondary: 13F55; 05E40.}

\end{minipage}
\end{center}

\section{1. Preliminaries}

Throughout, let $\mathfrak K$ be a field and let $S=\mathfrak K[x_1,\ldots,x_n]$ be the polynomial ring over $\mathfrak K$. For any square-free monomial ideal $I$ of $S$, let $G(I)$ be the set of minimal generator set of monomials, and let ${\rm sm}(I)$ be the set of square-free monomials. For the ideal $I$, there exist two related simplicial complexes, i.e., the nonface simplicial complex
$$\delta_{\mathcal{N}}(I)=:\{\,  F\in 2^{[\, n\, ]}\mid X_F\in {\rm sm}(S)\smallsetminus {\rm sm}(I) \,  \}$$
of $I$ and the facet simplicial complex
$$\delta_\mathcal{ F}(I)=:\langle\,  F\in 2^{[\, n\, ]}\mid X_F\in {\rm G}(I)\, \rangle$$
of the clutter $G(I)$. If they possess  a same $f$-vector, then the  ideal $I$ is called an {\it $f$-ideal}. For a graph $G$,  if its edge ideal
$I(G)=:\langle\,\{X_F\mid F\in E(G)\}\,\rangle$
is an $f$-ideal,  then $G$ is called an {\it $f$-graph}. Refer to \cite{2012Anwar,AMBZ2014,2014Anwar,GuoWuFideals2,2016Anwar} for further related studies.

For a simplicial complex $\D$ on the vertex set $[n]$, let $\mathcal F (\D)$ be the clutter of facets  in $\D$, let $\mathcal N(\D)$ be the set of all minimal nonfaces of $\D$. Let
$$\D^c=:\langle\{ F\mid F^c\in \mathcal F (\D)\}\rangle\overset{i.e.}{=}\langle \{[n]-G\mid G\in \mathcal F (\D)\}\rangle.$$
$\D$ is called an {\em $f$-simplicial complex} if the facet ideal $I(\D)$ is an $f$-simplicial complex. Note that in defining an $f$-graph $G$, $G$ is regarded as a simplicial complex of dimension no more than $1$, although we do have $I(G)=I_{{\rm Ind}(G)}$, where ${\rm Ind}(G)$ is the independence simplicial complex of the graph $G$. The definition of an $f$-simplicial complex seems to be reasonable with hindsight, due to the following two theorems on $f$-ideals.

\begin{thm}\label{homogeneous $f$-Ideal} {\rm (\cite[Theorem 2.3]{GuoWuLiu})} Let $S=K[x_1,  \ldots,  x_n]$,  and let $I$ be a  square-free monomial ideal of $S$ with the minimal generating set $G(I)$, where all monomials of $ G(I)$ have a same homogeneous degree $d$. Then $I$ is an $f$-ideal if and only if,  the set $ G(I)$ is an LU-set and,  $| G(I)| = \frac{1}{2}\binom{n}{d}$ holds true.
\end{thm}

Note that $ G(I)$ is said to be an LU-set if the set of all degree $d-1$ factors of elements of $G(I)$ has exactly $\binom{n}{d-1}$ elements, and the set of degree $d+1$ square-free monomials extended from elements of $G(I)$ has cardinality $\binom{n}{d+1}$.

Recall also the following recently discovered result:

\begin{thm} {\rm(\cite[Theorem 4.1]{BTuyl2018})}\label{Newton} $\D $ is an $f$-simplicial complex, if and only if $\D^c$ is an $f$-simplicial complex.
\end{thm}

Equivalently, a square-free monomial ideal $I$ of $S$ is an $f$-ideal if and only if the following  Newton complement dual ideal
$$\hat{I}=\langle x_1x_2\cdots x_n/u\mid u\in G(I)\rangle$$
of $I$ is an $f$-ideal.

It is clear that  Theorem \ref{Newton} follows easily from Theorem \ref{homogeneous $f$-Ideal} for a pure simplicial complex $\D$.

\section{2. Well-distributed $f$-simplicial complexes}

We begin with the following interesting example.

\begin{myeg} \label{f-idealNonShellable}
Consider the simplicial complex $\D$ whose facet set is
$$\{123,  125,  136,  145,  146,  234,  246,  256,  345,  356 \}.$$
By {\rm \cite[Example 7.7]{BW1996}}, it is not shellable.  It is direct to check that $\mathcal F (\D)$ is a LU-set and $|\mathcal {\rm G}(I)|=\frac{1}{2}\binom 6 3$. Thus the facet ideal $I=:I(\D)$ is an $f$-ideal, hence $\hat{I}=:I(\D^c)$ is also an $f$-ideal. Note that $\mathcal F (\D)\cap \mathcal F (\D^c)=\eset.$

By taking advantage of CoCoA, we get the primary decomposition of $I$ as follows:
$$I(\D)=\langle x_3, x_5, x_6\rangle\cap \langle x_2, x_4, x_6\rangle\cap\langle x_2, x_5, x_6\rangle\cap\langle x_1, x_4, x_6\rangle\cap\langle x_1, x_3, x_6\rangle\cap\langle x_3, x_4, x_5\rangle$$
$$\cap \langle x_1, x_4, x_5\rangle\cap\langle x_1, x_2, x_5\rangle\cap\langle x_2, x_3, x_4\rangle\cap\langle x_1, x_2, x_3\rangle,$$
hence $I$ is unmixed. Furthermore,  we get the following same $3$-linear resolution for both $I$ and $\hat{I}$ by CoCoA:
$$0\longrightarrow  S^6(-5)\longrightarrow S^{15}(-4)\longrightarrow  S^{10}(-3)\longrightarrow S.$$
Hence by $I(\D)=I_{(\D^c)^\vee}$, both ideals $I(\D)$ and $I(\D^c)$ are \CM by Eagon-Reiner theorem.
\end{myeg}

In order to seek more information of the example, we introduce the following concept:

\begin{mydef}
  Let $n=2d$, and let $\D$ be a pure $d-1$-dimensional simplicial complex over vertex set $[n]$. If $\mathcal{F}(\D)\cap \mathcal{F}(\D^c)=\eset$ \,holds, then $\D$ is called a  well-distributed  simplicial complex.
\end{mydef}

Let $\pi$ be a permutation on $[n]$, and let $\D$ be a pure simplicial complex with vertex set $[n]$. By Theorem \ref{homogeneous $f$-Ideal}, $\pi(\D)=:\{\pi(F)\mid F\in \D\}$ is an $f$-simplicial complex if $\D$ is. For a $(2d-1)$-dimensional pure $f$-simplicial complex $\D$, note that $\D$ is  well-distributed  \iff there exists no facet $F\in \D$ such that $[n]-F\in\mathcal F(\D)$, thus $\pi(\D)$ is still  well-distributed  if $\D$ is.

 We have the following observation.

\begin{mthm} \label{AnEquality}  For  a  well-distributed  $f$-simplicial complex $\D$, $\mathcal N(\D^c)=\mathcal F(\D)$ holds true.
\end{mthm}

\Proof Assume  ${\rm dim}(\D)=d-1$ and $V(\D)=[n]$.  By Theorem \ref{Newton} and \ref{homogeneous $f$-Ideal}, $\mathcal{F}(\D^c)$ is a L-set, thus  ${\rm sm}(S)_{d-1}\seq \D^c$. Then the condition $\mathcal{F}(\D)\cap \mathcal{F}(\D^c)=\eset$ implies $\mathcal F (\D)\seq \mathcal N(\D^c)$. Conversely, for any $U\seq [n]$ with $|U|>d$, since $\mathcal F (\D)$ is a U-set, we have $H\in \mathcal F (\D)$ such that $H\seq U$, hence $U\not\in \mathcal N(\D^c)$. Finally, the condition $\mathcal{F}(\D)\cap \mathcal{F}(\D^c)=\eset$ implies $\mathcal N(\D^c)=\mathcal F(\D)$.
\quad\quad\qed

\begin{mycor} \label{AlexanderDual}  For  a  well-distributed   $f$-simplicial complex $\D$ , $\D^\vee=\D$ holds true.
\end{mycor}

\Proof By Lemma \ref{AnEquality}, we have $\mathcal N(\D)=\mathcal F(\D^c)$.
 By definition of $\D^\vee$, we have $$\mathcal F (\D^\vee)=\{[n]-F\mid F\in \mathcal N(\D)\}=\{[n]-F\mid F\in \mathcal F(\D^c)\}=\mathcal F (\D).$$
Thus $\D^\vee=\D$.\quad\quad\qed

\vs{3mm}With the observations, Eagon-Reiner theorem has a stronger form for the small class of  well-distributed  $f$-simplicial complexes:

\begin{mthm} \label{$f$-idealCM}  For  a  well-distributed  $f$-simplicial complex $\D$, let $J=I(\D)$. Then \tfae:

$(1)$ The ideal $J$ is Cohen-Macaulay.

$(2)$ $J$ has linear resolution.

\end{mthm}

\Proof  By Lemmas \ref{AnEquality}, we have $I(\D)=I_{\D^c}$. Then apply Eagon-Reiner theorem to the equality $I(\D)=I_{(\D^c)^\vee}$, $\D^c$ is \CM \iff $I(\D)$ has linear resolution.  On the other hand, $\D^c$ is \CM \iff $I_{(\D^c)^\vee}$ is a \CM ideal. Since  $\D^c$ is also a  well-distributed  $f$-simplicial complex, we get
$I(\D)=\D^c=(\D^c)^\vee$ by Lemma \ref{AlexanderDual}. This completes the proof.\quad\quad \qed

\vs{3mm}
Now we can go back to answer the question posed after Example \ref{f-idealNonShellable}. Note that both $\D$ and $\D^c$ are  well-distributed  $f$-simplicial complexes, thus Theorem \ref{$f$-idealCM} applies. Hence  both $I(\D)$ and $I(\D^c)$ are \CM ideals. Note that though both $\D$ and $\D^c$ are \CM,  neither  is shellable.

Recall that in \cite{GuoWuLiu}, a  characterization (in fact,  a complete classification) of $f$-graphs was presented. Based on the classification, it was proved in \cite{GuoWuLiu} that all  $f$-graphs are connected, thus well-covered and \vd\, and hence, pure shellable. In particular,  all $f$-graphs are Cohen-Macaulay. In contrast, there exist a lot of nonpure $f$-simplicial complexes of dimension greater than $1$ (\cite{GuoWuLiu,GuoWuFideals2}).  Example \ref{f-idealNonShellable} gives more evidence showing that the higher dimensional $f$-simplicial complexes are a little bit complicated.

\section{3. Strong $f$-simplicial complexes}

Throughout this section, let $\D$ be a $(d-1)$-dimensional pure simplicial complex over vertex set $[n]$. Let $[n]_d=\{A\seq 2^{[n]}\mid |A|=d\}.$ Let $\D'=:\langle F\mid F\in [n]_d\smallsetminus \mathcal F(\D)\rangle$ be the homogeneous complement of $\D$.

Inspired by the concept of  well-distributed  $f$-simplicial complex, we now introduce the following definitions:

\begin{mydef} For  a pure simplicial complex $\Delta$, if both $\D$ and its homogeneous complement complex $\Delta'$ are $f$-simplicial complexes, then $\D$ is said to be {\bf strong.}  A square-free monomial ideal $I$ is called a  strong $f$-ideal, if its facet simplicial complex $I(\D)$ is strong.
\end{mydef}

Consider all $(d-1)$-dimensional pure simplicial complexes with  vertex set $[2d]$. $f$-ideals are abundant among the kind of complexes, but only a few are  well-distributed. For example, when $d=2$, there are totally $20$ such complexes, among them $12$ are $f$-simplicial complexes, but none is  well-distributed.

\begin{myeg} We list all the $1$-dimensional pure $f$-simplicial complexes with $|\mathcal F(\D)|=3$:
$$\langle 12,13,24  \rangle,\langle 12,13,34  \rangle,\langle  12,14,23 \rangle,\langle 12,14,34  \rangle,\langle 12,23,34  \rangle,\langle 12,24, 34  \rangle,$$
 $$\langle14,23,34\rangle,\langle 14,23,24   \rangle,  \langle 13,24,34  \rangle,\langle 13,23,24  \rangle, \langle 13,14,24  \rangle,\langle 13,14,23  \rangle.$$
Note that homogeneous complements of the $f$-simplicial complexes are all $f$-simplicial complexes, thus all are strong $f$-simplicial complexes. Note also that all the $f$-simplicial complexes are isomorphic, in the sense that they are all paths of lengths of three. Alternatively, there are totally two classes of the codimension 1 graphs of the simplicial complex, i.e., the star graph $K_{1,3}$ which corresponds to the eight non-$f$-simplicial complexes  and, the line graph $L_3$ which corresponds to the 12 $f$-simplicial complexes.

\end{myeg}

When $d=3$, there exist examples of unmixed $f$-simplicial complexes that is not strong, as the following example shows:

\begin{myeg} Let $\D$ be a simplicial complex on vertex set $[6]$ with facet set
$$\mathcal F(\D)=\{123,  124,  125,  126,  345,  346,  456,  356,  134,  256\}.$$
Then it is direct to check that $\mathcal F(\D)$ is a U-set and an L-set, thus $\D$ is an $f$-simplicial complex. The facet set of the homogeneous complement of $\D$ is
$$\{135,  136,  145,  146,  156,  234,  235,  236,  245,  246\}.$$
It is not an L-set since it does not cover $12$. Thus $\D'$ is not an $f$-simplicial complex, thus $\D$ is not strong.
\end{myeg}

Examples of strong $f$-simplicial complexes are abundant, see \cite{GuoWuFideals2} for all the constructions of several classes of pure $f$-simplicial complexes.

Next, we show that $^{'}$ and $^c$ commutes:
\begin{myprn}\label{'and c commutes} For any
pure simplicial complex $\D$, $(\D^c)'=(\D')^c$ holds true.
\end{myprn}
\Proof Let $[n]$ be the vertex set of $\D$ and assume ${\rm dim }(\D)=d$. Then for any $F\in 2^{[n]}$ with $|F|=n-d$, we have
$$F\in \mathcal F((\D')^c)\Llra [n]-F\in \D'\Llra [n]-F\not\in \mathcal F(\D)\Llra F\not\in \mathcal F(\D^c)\Llra F\in \mathcal F ((\D^c)').$$
Thus $(\D^c)'=(\D')^c$ holds true.
\quad\quad\qed

\vs{2mm}Clearly, the map $\vi: \,[n]_d\to [n]_{n-d}, \, F\mapsto [n]-F$ is bijective, thus the decomposition $[n]_d=\mathcal F({\D})\,\overset{\cdot}{\bigcup}\,\mathcal F({\D'})$ implies $\vi(\mathcal F({\D'}))=\mathcal F({(\D^c)'})$. Thus by Theorem \ref{Newton}, we have

\begin{prop} \label{N complement}  A simplicial complex $\D$ is strong \iff $\D^c$ is strong.
\end{prop}

We have the following observation:

\begin{prop} \label{h complement}  Let $\D$ be a strong $f$-simplicial complex. Then

$(1)$ $\mathcal N(\D)=\mathcal F(\D')$, hence we have
$I_{\D}= I(\D')$.

$(2)$ $ \D^\vee=(\D^c)'$.
\end{prop}

\Proof $(1)$ The equality $I_{\D}= I(\D')$ follows by definition of $\D'$ and the proof of Theorem \ref{AnEquality}.

$(2)$ Since $\mathcal N(\D)=\mathcal F(\D')$, by definition of $\D^\vee$, we have $$\mathcal F (\D^\vee)=\{[n]-F\mid F\in \mathcal N(\D)\}=\{[n]-F\mid F\in \mathcal F(\D')\}=\mathcal F ((\D')^c).$$
Thus by Proposition \ref{'and c commutes}, $\D^\vee=(\D^c)'$ holds true.
\quad\quad\qed

We have the following main result of this section:

\begin{thm} \label{$f$-idealCM1} For  a strong $f$-simplicial complex $\D$,  \tfae:

$(1)$ The ideal $I(\D')$ is Cohen-Macaulay.

$(2)$ $\D$ is a Cohen-Macaulay simplicial complex.

$(3)$ The ideal $I(\D^c)$ has linear resolution.

 \end{thm}

\Proof By Eagon-Reiner theorem, $I_\D$ is \CM \iff $\D$  is Cohen-Macaulay, the latter holds true  if and only if $I_{\D^\vee}$ has linear resolution. Then the result follows by applying Propositions \ref{h complement} and \ref{'and c commutes}. \quad\quad \qed

\begin{mycor} For  a strong $f$-simplicial complex $\D$,  \tfae:

$(1)$ The ideal $I(\D^c)$ is Cohen-Macaulay.

$(2)$ The ideal $I(\D')$ has linear resolution.

\end{mycor}


\begin{thebibliography}{gg}

\bibitem{2012Anwar} G.Q. Abbasi,  S. Ahmad,  I. Anwar and W.A. Baig. $F$-Ideals of degree 2. Algebra Colloq. $19 (2012)\,  921-926.$

\bibitem{AMBZ2014} I. Anwar,  H. Mahmood,  M. A. Binyamin and M. K. Zafar.  On the Characterization of $f$-Ideals,  Commun. Algebra $42(2014)$,  $3736-3741$.

\bibitem{BW1996} A. Bj\"orner,   M.L. Wachs.  Shellable nonpure complexes and posets. I,  Trans. Amer. Math. {\rm Soc}.,   $348(1996)\, \,  1299-1327.$

\bibitem{BTuyl2018} S. Budd,  A. Van Tuyl. Newton complementary duals of $f$-ideals. Canad. Math. Bull. $62:2(2019)\,  231-241.$
ArXiv: $1804.00686v1. $

\bibitem{GuoWuLiu} J. Guo,  T.S. Wu and Q. Liu.  $F$-ideals and $f$-graphs.   Comm. Algebra $45:8$ $(2017)\,  3207-3220$.

\bibitem{GuoWuFideals2} J. Guo,  T.S. Wu. On the  $(n,  d)^{th}$ $f$-ideals. J. Korean Math. {\rm Soc}. $52:4(2015)\,  685-697.$

\bibitem{Herzog and Hibi} J. Herzog and  T. Hibi. {\it Monomial Ideals}.  GTM $260$ London: Springer-Verlag London Limited,  $2011$.

\bibitem{2016Anwar} H. Mahmood,  I. Anwar,  M.A. Banyamin and S. Yasmeen. On the connectedness of $f$-simplicial complexes. J. Algebra Appl. $15:6(2016)\,1750017,  9 pp.$  doi: $10.1142/S0219498817500177.$

\bibitem{2014Anwar} H. Mahmood,  I. Anwar and M.K. Zafar. Construction of Cohen-Macaulay f-Graphs.
J. Algebra Appl. $13:6(2014)\, 14500121,  7 pp.$


\bibitem{Moradi2} R. Rahmati-Asghar and S. Moradi. On the Stanley-Reisner ideal of an expanded simplicial complex. Manuscripta Math. $150:3-4 (2016)\, \,  533-545.$

\end{thebibliography}
\end{document}